\documentclass{article}    
\usepackage{amstext}
\def\qed{$\rlap{$\sqcap$}\sqcup$}
\usepackage{latexsym}

\begin{document}           

{\ }\\
\begin{center}
{\huge {\bf Partial derivatives of a generic subspace of a vector space of forms: quotients of level algebras of arbitrary type}} \\ [.250in]
{\large FABRIZIO ZANELLO\\
Dipartimento di Matematica, Universit\`a di Genova, Genova, Italy.\\E-mail: zanello@dima.unige.it}
\end{center}

{\ }\\
\\
ABSTRACT. Given a vector space $V$ of homogeneous polynomials of the same degree over an infinite field, consider a generic subspace $W$ of $V$. The main result of this paper is a lower-bound (in general sharp) for the dimensions of the spaces spanned in each degree by the partial derivatives of the forms generating $W$, in terms of the dimensions of the spaces spanned by the partial derivatives of the forms generating the original space $V$.\\
Rephrasing our result in the language of commutative algebra (where this result finds its most important applications), we have: let $A$ be a type $t$ artinian level algebra with $h$-vector $h=(1,h_1,h_2,...,h_e)$, and let, for $c=1,2,...,t-1$, $H^{c,gen}=(1,H_1^{c,gen},H_2^{c,gen},...,H_e^{c,gen})$ be the $h$-vector of the generic type $c$ level quotient of $A$ having the same socle degree $e$. Then we supply a lower-bound (in general sharp) for the $h$-vector $H^{c,gen}$. Explicitly, we will show that, for any $u\in \lbrace 1,...,e\rbrace $, $$H_u^{c,gen}\geq {1\over t^2-1}\left((t-c)h_{e-u}+(ct-1)h_u\right).$$
This result generalizes a recent theorem of Iarrobino (which treats the case $t=2$).\\
Finally, we begin to obtain, as a consequence, some structure theorems for level $h$-vectors of type bigger than 2, which is, at this time, a very little explored topic.\\
\\
Math. Subject Class.: 13E10 (Primary); 13H10 (Secondary).\\
Key words and phrases: Artinian algebra; level algebra; $h$-vector; generic quotient; dimension; partial derivatives.\\
\\
\section{Introduction and preliminary results}
\indent
{\large

Let $S$ be a polynomial ring in $r$ variables over an infinite field, and denote by $S_i$ the subspace of degree $i$ homogeneous elements (forms) of $S$. Let us consider a vector subspace $V$ of $S_e$, generated by $t$ forms, and let $h=(1,h_1,...,h_{e-1},h_e=t)$ be the vector of the dimensions of the spaces spanned by the partial derivatives of the generators of $V$ in each degree. That is, for $i=0,1,...,e$, $h_{e-i}$ is the dimension of the subspace of $S_{e-i}$ spanned by the order $i$ partial derivatives of the elements of $V$. Let $W\subseteq V$ be a vector subspace generated by $c$ {\it generic} (this term will be defined later) forms of $V$. The main result of this paper is a lower-bound (which is in general sharp) for the vector $H^{c,gen}=(1,H_1^{c,gen},H_2^{c,gen},...,H_e^{c,gen})$ of the dimensions of the spaces spanned by the partial derivatives of the generators of $W$. This lower-bound will be expressed in terms of the entries of the vector $h$ above.\\\indent
The most important application of this result is to commutative algebra, in particular to the study of the $h$-vectors of level algebras, as we will see later. We remark here that there are also significant applications of the study of the behavior of the partial derivatives of multivariate polynomials to other areas of algebra (e.g., invariant theory), to geometry, combinatorics, and even complexity theory (see, e.g., the article $[NW]$).\\
\\\indent
Before rephrasing the result mentioned above in the language of commutative algebra, let us introduce the definitions we will need and fix the setting we will be working in throughout this paper. We consider standard graded artinian algebras $A=R/I$, where $R=k[x_1,...,x_r]$, $I$ is a homogeneous ideal of $R$ and the $x_i$'s all have degree 1. We assume for simplicity that $k$ is a field of characteristic zero, even if our results also hold when $k$ is an infinite field of characteristic $p$ (see Remark 2.15).\\\indent
The {\it $h$-vector} of $A$ is $h(A)=h=(h_0,h_1,...,h_e)$, where $h_i=\dim_k A_i$ and $e$ is the last index such that $\dim_k A_e>0$. Since we may suppose that $I$ does not contain non-zero forms of degree 1, $r=h_1$ is defined to be the {\it codimension} of $A$.\\\indent 
The {\it socle} of $A$ is the annihilator of the maximal homogeneous ideal $\overline{m}=(\overline{x_1},...,\overline{x_r})\subseteq A$, namely $soc(A)=\lbrace a\in A {\ } \mid {\ } a\overline{m}=0\rbrace $. Since $soc(A)$ is a homogeneous ideal, we define the {\it socle-vector} of $A$ as $s(A)=s=(s_0,s_1,...,s_e)$, where $s_i=\dim_k soc(A)_i$. Notice that $h_0=1$, $s_0=0$ and $s_e=h_e>0$. The integer $e$ is called the {\it socle degree} of $A$ (or of $h$). The {\it type} of the socle-vector $s$ (or of the algebra $A$) is type$(s)=\sum_{i=0}^es_i$.\\\indent 
If $s=(0,0,...,0,s_e=t)$, we say that the algebra $A$ is {\it level} (of type $t$). In particular, if $t=1$, $A$ is {\it Gorenstein}. With a slight abuse of notation, we will sometimes refer to an $h$-vector as Gorenstein (or level) if it is the $h$-vector of a Gorenstein (or level) algebra.\\\indent
We now recall the main facts of the theory of {\it inverse systems}, or {\it Macaulay duality}, which we will use in the sequel. Inverse systems are the key instrument which links our main theorem to commutative algebra. (For a complete introduction to this theory, we refer the reader to $[Ge]$ and $[IK]$.)\\
\\\indent
{\bf Macaulay duality.} Let $S=k[y_1,...,y_r]$, and consider $S$ as a graded $R$-module where the action of $x_i$ on $S$ is partial differentiation with respect to $y_i$.\\\indent 
There is a one-to-one correspondence between artinian algebras $R/I$ and finitely generated $R$-submodules $M$ of $S$, where $I=Ann(M)$ is the annihilator of $M$ in $R$ and, conversely, $M=I^{-1}$ is the $R$-submodule of $S$ which is annihilated by $I$ (cf. $[Ge]$, Remark 1), p. 17).\\\indent 
If $R/I$ has socle-vector $s$, then $M$ is minimally generated by $s_i$ elements of degree $i$, for $i=1,2,...,e$, and the $h$-vector of $R/I$ is given by the number of linearly independent derivatives in each degree obtained by differentiating the generators of $M$ (cf. $[Ge]$, Remark 2), p. 17).\\\indent 
In particular, level algebras of type $t$ and socle degree $e$ correspond to $R$-submodules of $S$ minimally generated by $t$ elements of degree $e$.\\ 
\\\indent
The present work generalizes the main result of Iarrobino's paper $[Ia]$, which determines a lower-bound (that is in general sharp, as we have noticed in $[Za]$) for the $h$-vector of the generic Gorenstein quotient of a type 2 level $k$-algebra $A$ of the same socle degree, say $e$. In fact, it is easy to see, by inverse systems, that the generic Gorenstein quotients of $A$ of socle degree $e$ all have the same $h$-vector. (Here, in the type 2 context, by {\it generic} Gorenstein quotients we mean, as in Iarrobino's work, those parameterized by the points of a non-empty open subset of ${\bf P}^1(k)$, given two generators for the inverse system module corresponding to $A$.)\\\indent
Furthermore, the $h$-vector of those generic Gorenstein quotients is the (entry by entry) largest among the $h$-vectors of {\it all} the Gorenstein quotients of $A$ of socle degree $e$. The key portion of Iarrobino's work was obtaining a \lq \lq good" lower-bound for this generic $h$-vector.\\\indent 
Since, by a theorem of Geramita {\it et al.} ($[GHMS]$) that we will state later (Theorem 1.5), the reverse of the difference between the $h$-vector of a level algebra of type $t$ and the $h$-vector of any of its type $t-1$ quotients having the same socle degree is an $O$-sequence, Iarrobino was able to deduce from his theorem some restrictions on the possible $h$-vectors of level algebras of type $t=2$ (we also extensively exploited $[GHMS]$'s result in $[Za]$).\\\indent
Although level $h$-vectors have been studied a lot ($[GHMS]$'s bibliography lists all the main papers on the subject up to 2003), at this date there is no complete characterization of these $h$-vectors, even in the first non-Gorenstein case: type 2 and codimension 3 (see $[Ia2]$ for codimension 2). For $t>2$, level $h$-vectors are even less known (see, e.g., the list of problems in $[GHMS]$.)\\\indent
Given an artinian level algebra $A$ of arbitrary type $t$, the most important result of this paper is a lower-bound (in general sharp) for the $h$-vector of the generic type $c$ ($c=1,2,...,t-1$) level quotient of $A$ having the same socle degree (Main Theorem (2.9)).\\\indent
We will conclude by exhibiting some applications, which show that our result is a useful tool in starting to study effectively the still wide open problem of classifying level $h$-vectors of arbitrary type.\\\indent
On the other hand, of course, in this paper we will not yet be able to say for type $t>2$ as much as can be said for type 2, since Gorenstein $h$-vectors are much better known than type $t-1$ level $h$-vectors when $t>2$. Therefore, we will be able to exploit the main result of this work in the study of type $t$ level $h$-vectors better when we know more about type $t-1$ level $h$-vectors.\\
\\\indent
We now state the main result of Iarrobino's paper $[Ia]$ (which we rephrase in our notation):\\
\\\indent
{\bf Theorem 1.1} ($[Ia]$). {\it Let $A$ be a level algebra of type 2 and socle degree $e$ having $h$-vector $h=(1,h_1,h_2...,h_e)$, and let $H^{1,gen}=(1,H_1^{1,gen},H_2^{1,gen},...,H_e^{1,gen})$ be the $h$-vector of the generic Gorenstein quotient of $A$ having socle degree $e$. Furthermore, let $u\in \lbrace 1,2,...,e\rbrace $, and define $i=e-u$. If, for some integer $\delta_u\geq 0$, we have $h_i\geq 2h_u-2-3\delta_u$, then $H_u^{1,gen}\geq h_u-\delta_u$.}\\
\\\indent
{\bf Proof.} See $[Ia]$, Theorem 2.4.{\ }{\ }\qed \\
\\\indent
{\bf Remark 1.2.} Notice that Iarrobino's theorem more simply states that, with the notation above, \begin{equation}\label{iii}H_u^{1,gen}\geq {1\over 3}(h_u+h_i).\end{equation}\indent
This equivalent form (which gets rid of $\delta_u$ and will be consistent with our following formulas) can be easily obtained by observing that the smallest integer $\delta_u$ satisfying the inequality of the statement is $\delta_u=\left\lceil {2h_u-h_i-2\over 3}\right\rceil $, where as usual $\left\lceil \alpha \right\rceil $ indicates the least integer greater than or equal to $\alpha $ and $\lfloor \alpha \rfloor $ indicates the largest integer less than or equal to $\alpha $. Therefore the conclusion of the theorem is $H_u^{1,gen}\geq h_u- \left\lceil {2h_u-h_i-2\over 3}\right\rceil =h_u- \left\lfloor {2h_u-h_i\over 3}\right\rfloor =\left\lceil {1\over 3}(h_u+h_i)\right\rceil $, whence we obtain (\ref{iii}), since $H_u^{1,gen}$ is an integer.\\
\\\indent
Using the above notation, our generalization of Iarrobino's result (in the form of inequality (\ref{iii})), which we now state in the language of commutative algebra, is the following (we will prove it in the next section):\\
\\\indent
{\bf Theorem 2.9 (Main Theorem).} {\it For every integer $c$, $1\leq c\leq t-1$, we have: $$H_u^{c,gen}\geq {1\over t^2-1}\left((t-c)h_i+(ct-1)h_u\right).$$}
\\\indent
{\bf Definition-Remark 1.3.} Let $n$ and $i$ be positive integers. The {\it i-binomial expansion of n} is $$n_{(i)}={n_i\choose i}+{n_{i-1}\choose i-1}+...+{n_j\choose j},$$ where $n_i>n_{i-1}>...>n_j\geq j\geq 1$. Under these hypotheses, the $i$-binomial expansion of $n$ is unique (e.g., see $[BH]$, Lemma 4.2.6).\\\indent 
Furthermore, define $$n^{<i>}={n_i+1\choose i+1}+{n_{i-1}+1\choose i-1+1}+...+{n_j+1\choose j+1}.$$
\\\indent 
{\bf Theorem 1.4} (Macaulay). {\it Let $h=(h_i)_{i\geq 0}$ be a sequence of non-negative integers, such that $h_0=1$, $h_1=r$ and $h_i=0$ for $i>e$. Then $h$ is the $h$-vector of some standard graded artinian algebra if and only if, for every $d$, $1\leq d\leq e-1$, $$h_{d+1}\leq h_d^{<d>}.$$}
\\\indent 
{\bf Proof.} See $[BH]$, Theorem 4.2.10. (This theorem holds, with appropriate modifications, for any standard graded algebra, not necessarily artinian.){\ }{\ }\qed \\
\\\indent
A sequence of non-negative integers which satisfies the growth condition of Macaulay's theorem is called an {\it $O$-sequence}.\\\indent
The following is the fundamental result of Geramita {\it et al.} that we mentioned earlier:\\
\\\indent
{\bf Theorem 1.5} ($[GHMS]$). {\it Let $F_1,...,F_t\in S$ be any linearly independent forms of degree $e$, and let $M=<F_1,...,F_t>$ and $N=<F_1,...,F_{t-1}>$ be two inverse system submodules of $S$. Then the reverse of the difference between the $h$-vectors of $A=R/Ann(M)$ and $R/Ann(N)$ is an $O$-sequence (which is the $h$-vector of a quotient of $A$).}\\
\\\indent
{\bf Proof.} See $[GHMS]$, Lemma 2.8 and Theorem 2.10.{\ }{\ }\qed \\
\\
\section{The main result}
\indent

Let us consider a level algebra $A=R/I$ having type $t$, socle degree $e$ and $h$-vector $h=(1,h_1,h_2...,h_e)$, and let $M=<F_1,F_2,...,F_t>\subseteq S$ be the inverse system module which is annihilated by $I$. It can be easily seen that, similarly to what we already observed in the Introduction for $t=2$, for any $t\geq 2$ and $c\in \lbrace 1,2,...,t-1\rbrace $ the {\it generic} level quotients of $A$ of type $c$ and socle degree $e$ (that is, the quotients parameterized by the points of a non-empty open subset of $({\bf P}^{t-1}(k))^c$, given the generators $F_1$, ..., $F_t$ of $M$) all have the same $h$-vector, which is the (entry by entry) maximal $h$-vector among those of {\it all} the type $c$ level quotients of $A$ having socle degree $e$. Our purpose is to determine a \lq \lq good" lower-bound for this generic $h$-vector, $H^{c,gen}=(1,H_1^{c,gen},H_2^{c,gen},...,H_e^{c,gen})$.\\\indent
Recall that, by inverse systems (see the Macaulay duality section in the Introduction), the integer $H_i^{c,gen}$ is the dimension of the space spanned by the $(e-i)$-th partial derivatives of the forms generating a generic $c$-dimensional subspace of $M$ in degree $e$ .\\
\\\indent
{\bf Definition-Remark 2.1.} As above, let $i=e-u$. For given forms $G_1,G_2,...,G_p\in S$ of degree $e$, define, for any $u$ such that $1\leq u\leq e-1$, $d_u(G_1,G_2,...,G_p)$ as the dimension (over $k$) of the intersection of the spaces spanned by the $i$-th partial derivatives of the $G_j$'s (which is clearly a $k$-vector subspace of $S_u$).\\\indent
Given $M\subseteq S$ and a minimal set of generators $F_1,...,F_t$ of $M$, define $$\Sigma_u(F_1,...,F_t)=\sum_{j_1<j_2}d_u(F_{j_1},F_{j_2})-\sum_{j_1<j_2<j_3}d_u(F_{j_1},F_{j_2},F_{j_3})+...+(-1)^td_u(F_1,...,F_t).$$\indent
(If the choice of the $F_j$'s is clear from the context, we will simply write $\Sigma_u$ instead of $\Sigma_u(F_1,...,F_t)$.)\\
\\\indent
{\bf Lemma 2.2.} {\it If $F_1,...,F_t$ are generic forms of $M=I^{-1}$ generating $M$, then: \begin{equation}\label{21}\Sigma_u=tH_u^{1,gen}-h_u.\end{equation}}
\\\indent
{\bf Proof.} Notice that, since the $F_j$'s are being chosen generically, the dimension of the quotient space given by the intersection of the spaces spanned by the $i$-th derivatives of, say, $q$ of these $t$ forms modulo (its intersection with) the space spanned by the derivatives of the other $t-q$ forms, is independent of the subset of $q$ forms that we choose. (We will call this dimension $D_u(q)$, for $q=1,...,t$. For instance, if $t=3$, the dimension $D_u(2)$ of the intersection of the spaces spanned by the $i$-th derivatives of $F_1$ and $F_2$ modulo the space spanned by the derivatives of $F_3$ remains constant when we permute the indices 1, 2 and 3.)\\\indent
The formula of the statement now immediately follows from a straightforward application of the exclusion-inclusion principle (see the next example).{\ }{\ }\qed \\
\\\indent
{\bf Example 2.3.} Let $t=3$. Then, by the exclusion-inclusion principle, $H_u^{1,gen}$ is equal to the following sum: the dimension $D_u(1)$ of the space spanned by the $i$-th derivatives of one generic form, say $F_1$, modulo the space of the derivatives of the other two forms ($F_2$ and $F_3$), plus the dimension $D_u(2)$ of the intersection of the spaces spanned by the derivatives of $F_1$ and $F_2$ modulo the space of the derivatives of $F_3$, plus one more time the dimension $D_u(2)$ (computed as above after switching $F_2$ and $F_3$), plus the dimension $D_u(3)$ of the intersection of the spaces spanned by the derivatives of $F_1$, $F_2$ and $F_3$. Thus, $H_u^{1,gen}=D_u(1)+2D_u(2)+D_u(3)$.\\\indent
Again by the exclusion-inclusion principle, one similarly checks that $h_u=3D_u(1)+3D_u(2)+D_u(3)$. Likewise, we have $\Sigma_u=\Sigma_u(F_1,F_2,F_3)=d_u(F_1,F_2)+d_u(F_1,F_3)+d_u(F_2,F_3)-d_u(F_1,F_2,F_3)=3(D_u(2)+D_u(3))-D_u(3)=3D_u(2)+2D_u(3)$.\\\indent
Hence it follows that, for $t=3$, $\Sigma_u=3H_u^{1,gen}-h_u$, as stated in the previous lemma.\\
\\\indent
Let us now supply a lower-bound for the $u$-th entry of the $h$-vector of the generic Gorenstein quotient of $A$ of socle degree $e$.\\
\\\indent
{\bf Theorem 2.4.} {\it For generic forms $F_1,F_2,...,F_t$ of $M$ generating $M$, we have \begin{equation}\label{eee}H_u^{1,gen}\geq h_i-\Sigma_u.\end{equation}}
\\\indent
{\bf Proof.} We proceed by induction on $t$, since we know that the result is true for $t=2$ ((\ref{eee}) is exactly the first inequality of the statement of $[Ia]$, Theorem 2.2 when $t=2$). Suppose that the theorem holds for $t-1$, where $t\geq 3$. Since generic Gorenstein quotients of $A=R/Ann(M)$ are clearly generic Gorenstein quotients of generic type $t-1$ level quotients of $A$, by the induction hypothesis we have \begin{equation}\label{1}H_u^{1,gen}\geq H_i^{t-1,gen}-\Sigma_u^{t-1},\end{equation} where by $\Sigma_u^{t-1}$ we indicate the integer $\Sigma_u$ referred to the generic type $t-1$ subalgebra.\\\indent
It is easy to see, by definition, that \begin{equation}\label{2}\Sigma_u^{t-1}=\Sigma_u-H_u^{1,gen}+D_u(1).\end{equation}\indent
It immediately follows from the definition of $D_i(1)$ that \begin{equation}\label{3}H_i^{t-1,gen}=h_i-D_i(1).\end{equation}\indent
Therefore, by (\ref{1}), (\ref{2}) and (\ref{3}) we have $$H_u^{1,gen}\geq h_i-\Sigma_u+H_u^{1,gen}-D_u(1)-D_i(1).$$\indent
Hence, in order to prove (\ref{eee}), it suffices to show that $$H_u^{1,gen}\geq D_u(1)+D_i(1).$$\indent
We will again proceed by induction on $t$, since for $t=2$ the result follows immediately from the first inequality of $[Ia]$, Theorem 2.2 (which coincides with (\ref{eee}) when $t=2$) and the symmetry of Gorenstein $h$-vectors. The dimensions $D_u^{t-1}(1)$ and $D_i^{t-1}(1)$, defined in the obvious way for the generic type $t-1$ level subalgebra, are clearly greater than or equal to, respectively, $D_u(1)$ and $D_i(1)$, for in computing the latter dimensions we must also take into account the intersections with the derivatives of one more form.\\\indent
Therefore, by induction and what we have just observed, we obtain $$H_u^{1,gen}\geq D_u^{t-1}(1) +D_i^{t-1}(1)\geq D_u(1)+D_i(1),$$ as we wanted to show. This proves the theorem.{\ }{\ }\qed \\
\\\indent
{\bf Corollary 2.5.} {\it $$H_u^{1,gen}\geq {1\over t+1}(h_i+h_u).$$}
\\\indent
{\bf Proof.} If we consider the $F_j$'s generic inside $M$, from Theorem 2.4 and equation (\ref{21}), we have $H_u^{1,gen}\geq h_i-\Sigma_u=h_i-tH_u^{1,gen}+h_u$, and the result immediately follows.{\ }{\ }\qed \\
\\\indent
The following corollary generalizes $[Ia]$, Corollary 2.5 to arbitrary $t$:\\
\\\indent
{\bf Corollary 2.6.} {\it Let $A$ be as above, $h(A)=(1,r,h_2...,h_{e-1},t)$. If $h_{e-1}\geq t(r-1)$, then there exists a Gorenstein quotient of $A$ of socle degree $e$ having codimension $r$.}\\
\\\indent
{\bf Proof.} By Corollary 2.5, for $h_{e-1}\geq t(r-1)$ we have $$H_1^{1,gen}\geq {1\over t+1}(h_{e-1}+r)\geq {1\over t+1}(t(r-1)+r)=r-{t\over t+1}>r-1,$$ whence $H_1^{1,gen}=r$.{\ }{\ }\qed \\
\\\indent
We now need the following combinatorial formula. We thank Professor A. Conca for providing us with the elegant proof below.\\
\\\indent
{\bf Lemma 2.7.} {\it Let $t$ and $j$ be integers such that $t\geq j\geq 2$. Then:
$$\sum_{h=2}^j(-1)^h{t\choose h}{t-h\choose j-h}=(j-1){t\choose j}.$$}
\\\indent
{\bf Proof.} An easy computation shows that $(-1)^0{t\choose 0}{t-0\choose j-0}+(-1)^1{t\choose 1}{t-1\choose j-1}=-(j-1){t\choose j}.$ Hence it suffices to show that \begin{equation}\label{cc}\sum_{h=0}^j(-1)^h{t\choose h}{t-h\choose j-h}=0.\end{equation}\indent
Let us fix $j$, and consider the summands of the l.h.s. of (\ref{cc}), for every $h$, as functions of $t$. They are equal to $$(-1)^{h}\cdot {t(t-1)\cdot \cdot \cdot (t-h+1)\over h!}\cdot {(t-h)(t-h-1)\cdot \cdot \cdot (t-j+1)\over (j-h)!}.$$\indent
These are clearly all functions of degree $h+(j-h)=j$, and each annihilates at $t=0,1,...,j-1$. Hence the l.h.s. of (\ref{cc}) is a function of $t$, say $f_j(t)$, of degree at most $j$ with the $j$ zeroes $t=0$, $t=1$, ..., $t=j-1$. In order to prove that $f_j(t)$ is identically equal to 0, it is enough to show that it has one more zero. Let us compute $f_j(j)$.
We have $$f_j(j)=\sum_{h=0}^j(-1)^h{j\choose h}{j-h\choose j-h}=\sum_{h=0}^j(-1)^h{j\choose h}=(1-1)^j=0.$$\indent
This proves (\ref{cc}) and the lemma.{\ }{\ }\qed \\
\\\indent
Using the same notation as above, notice that, for the type $t-1$ generic level subalgebra of $A$ of socle degree $e$, equation (\ref{21}) becomes \begin{equation}\label{t-1}\Sigma_u^{t-1}=(t-1)H_u^{1,gen}-H_u^{t-1,gen}.\end{equation}
\\\indent
{\bf Lemma 2.8.} {\it Let $t>2$. For generic forms $F_1,F_2,...,F_t$ of $M$ generating $M$, we have $$\Sigma_u={t-1\over t-2}\Sigma_u^{t-1}+c_1,$$ for some real number $c_1\geq 0$.}\\
\\\indent
{\bf Proof.} We begin by showing the following:\\\indent
{\it Claim.} For every $u$, we have \begin{equation}\label{ss}\Sigma_u={t\choose 2}D_u(2)+2{t\choose 3}D_u(3)+...+(t-2){t\choose t-1}D_u(t-1)+(t-1){t\choose t}D_u(t).\end{equation}\indent
{\it Proof of claim.} By definition of $\Sigma_u$, it is easy to check that the number of times $D_u(j)$ must be counted in computing $\Sigma_u$ is equal, for every $j=2,...,t$, to \begin{equation}\label{du}\sum_{h=2}^j(-1)^h{t\choose h}{t-h\choose j-h}.\end{equation}\indent
Thus, (\ref{ss}) follows from Lemma 2.7. This proves the claim.\\\indent
{\it Claim.} $$\Sigma_u^{t-1}=\left({t\choose 2}-{t-1\choose 1}\right)D_u(2)+\left(2{t\choose 3}-{t-1\choose 2}\right)D_u(3)+...$$$$+\left((t-2){t\choose t-1}-{t-1\choose t-2}\right)D_u(t-1)+\left((t-1){t\choose t}-{t-1\choose t-1}\right)D_u(t).$$\indent
{\it Proof of claim.} Let us suppose that, in computing $H^{t-1,gen}$, we consider the generic forms $F_1,F_2,...,F_{t-1}$ of $M$, and that we again consider these forms as the generators of the $t-1$ cyclic modules each giving the generic Gorenstein $h$-vector $H^{1,gen}$. The number of times we have to consider, in computing $\Sigma_u^{t-1}$, a quotient space where the derivatives of exactly $j$ of all the $t$ forms \lq \lq appear" (i.e. the coefficient of $D_u(j)$ in $\Sigma_u^{t-1}$) is clearly the number of times such a quotient space appears in $\Sigma_u$ minus the number of times one of these $j$ forms is $F_t$, i.e. the generic form we have \lq \lq left out" in computing $H_u^{t-1,gen}$. And this number is exactly ${t-1\choose j-1}$. This proves the claim.\\\indent
It is easy to check that, for every $j$, $2\leq j\leq t-1$, the ratio of the coefficients of $D_u(j)$ in $\Sigma_u$ and $\Sigma_u^{t-1}$ is greater than the ratio of the corresponding coefficients of $D_u(j+1)$, i.e. that $${(j-1){t\choose j}\over (j-1){t\choose j}-{t-1\choose j-1}}>{j{t\choose j+1}\over j{t\choose j+1}-{t-1\choose j}}.$$\indent
The result now immediately follows, since the ratio of the coefficients of $D_u(t)$ is ${t-1\over t-2}$.{\ }{\ }\qed \\
\\\indent
Now we are ready for the main result of this paper:\\
\\\indent
{\bf Theorem 2.9 (Main Theorem).} {\it For every integer $c$, $1\leq c\leq t-1$, we have: $$H_u^{c,gen}\geq {1\over t^2-1}\left((t-c)h_i+(ct-1)h_u\right).$$\indent
In particular, for $c=t-1$, $$H_u^{t-1,gen}\geq {1\over t^2-1}(h_i+(t^2-t-1)h_u).$$}
\\\indent
{\bf Proof.} Let us first show the theorem for $c=t-1$. We can assume that $t>2$, since the case $t=2$ is Iarrobino's theorem (rephrased as in (\ref{iii})). For the sake of brevity, let $c_0={t-1\over t-2}$.\\\indent
From Lemma 2.8 and equation (\ref{21}), we get $$\Sigma_u^{t-1}={1\over c_0}(\Sigma_u-c_1)={1\over c_0}(tH_u^{1,gen}-h_u-c_1).$$\indent Hence, (\ref{t-1}) becomes $$H_u^{t-1,gen}=(t-1)H_u^{1,gen}-{t\over c_0}H_u^{1,gen}+{1\over c_0}h_u+{c_1\over c_0}.$$\indent
Getting rid of ${c_1\over c_0}$, which is non-negative, and using Corollary 2.5, we obtain 
$$H_u^{t-1,gen}\geq {c_0t-c_0-t\over c_0}\cdot {1\over t+1}(h_i+h_u)+{1\over c_0}h_u,$$
which an easy computation shows is equal to ${1\over t^2-1}(h_i+(t^2-t-1)h_u)$, as we desired. This proves the theorem for $c=t-1$.\\\indent
Now let $c\leq t-2$. By induction, let us suppose that the theorem holds for some index $c+1$, $1\leq c\leq t-2$, i.e. that $$H_u^{c+1,gen}\geq {1\over t^2-1}\left((t-(c+1))h_i+(t(c+1)-1)h_u\right).$$\indent
We want to show the inequality of the statement for $H_u^{c,gen}$. Since generic level quotients of $A$ of type $c$ can be seen as generic level quotients of type $c$ of generic level quotients of $A$ of type $c+1$, we have
$$H_u^{c,gen}\geq {1\over (c+1)^2-1}(H_i^{c+1,gen}+(c(c+1)-1)H_u^{c+1,gen})\geq $$$${1\over (c+1)^2-1}\Bigg({1\over t^2-1}\left((t-(c+1))h_u+(t(c+1)-1\right)h_i)$$$$+(c(c+1)-1){1\over t^2-1}\left(\left(t-(c+1))h_i+(t(c+1)-1\right)h_u\right)\Bigg).$$\indent
A standard computation shows that the r.h.s. of the inequality above is equal to ${1\over t^2-1}((t-c)h_i+(ct-1)h_u),$ and this concludes the induction process. Since we have already shown the theorem for $c=t-1$, the proof is complete.{\ }{\ }\qed \\
\\\indent
Notice that Corollary 2.5 is the particular case $c=1$ of the Main Theorem (2.9).\\\indent
From the Main Theorem we immediately obtain the following:\\
\\\indent
{\bf Corollary 2.10.} {\it Let $A$ be as above, $h(A)=(1,r,h_2...,h_{e-1},t)$. If $h_{e-1}\geq -t^2+rt+2$, then there exists a level quotient of $A$ of type $t-1$ and socle degree $e$ having codimension $r$.}\\
\\\indent
{\bf Proof.} By the Main Theorem (2.9), for $h_{e-1}\geq -t^2+rt+2$ we have $$H_1^{t-1,gen}\geq {1\over t^2-1}(h_{e-1}+(t^2-t-1)r)\geq {1\over t^2-1}(-t^2+rt+2+(t^2-t-1)r)=r-{t^2-2\over t^2-1}>r-1,$$ and therefore $H_1^{1,gen}=r$.{\ }{\ }\qed \\
\\\indent
As we claimed in the Introduction, in general the Main Theorem (2.9) is sharp, that is there exist level algebras whose generic level subalgebras of type $c$ and socle degree $e$ have an $h$-vector which coincides with the lower-bound supplied by the theorem. We will see this in the next example. We will omit the details, since it extends {\it mutatis mutandis} $[Za]$, Remark 2.6 to any $t>2$ and $c>1$.\\
\\\indent
{\bf Example 2.11.} Let $A$ be the level algebra of type $t$ and codimension $r=(t+1)p$ associated to the inverse system module $M=<F_1,...,F_t>$, where $F_j=y_{jp+1}y_1^{e-1}+y_{jp+2}y_2^{e-1}+...+y_{(j+1)p}y_p^{e-1}.$\\\indent
It is easy to see that the $h$-vector of $A$ is $h=(1,(t+1)p,(t+1)p,...,(t+1)p,t)$. Standard considerations also show that all the level quotients of $A$ of type $c$ and the same socle degree (not only those generic) have $h$-vector $(1,(c+1)p,(c+1)p,...,(c+1)p,c)$, for every $c=1,2,...,t-1$.\\\indent
Since $$(c+1)p={1\over t^2-1}((t-c)(t+1)p+(tc-1)(t+1)p),$$ we have that the lower-bound of the Main Theorem (2.9) is sharp for these algebras $A$.\\
\\\indent
Although the Main Theorem (2.9), as we have seen, in general cannot be improved, there are cases where our lower-bound is not sharp. In particular, if we want a lower-bound for the $h$-vector of the generic type $c$ level subalgebra $B$ of $A$, sometimes we obtain a better result by considering $B$ as a subalgebra of a generic type $c+1$ level subalgebra of $A$, rather than directly as a subalgebra of $A$. The next example clarifies this concept.\\
\\\indent
{\bf Example 2.12.} Let $h=(1,3,5,7,7,5,3)$ be the $h$-vector of a level algebra $A$ of type 3. (Notice that such an $h$ exists: it suffices to truncate a Gorenstein algebra having $h$-vector $(1,3,5,7,7,5,3,1)$.)\\\indent
By the Main Theorem (2.9), a generic level subalgebra $B$ of $A$ of type 2 and the same socle degree has $h$-vector $h^{'}\geq (1,3,4,6,5,4,2)$. By $[GHMS]$'s Theorem 1.5 and Macaulay's Theorem 1.4, we actually have $h^{'}\geq (1,3,5,6,6,4,2)$ (since the reverse of the difference between $h$ and $h^{'}$ has to be an $O$-sequence, and (1,3,4,6) does not satisfy Macaulay's theorem).\\\indent
By applying the Main Theorem (2.9) to $h$, we have that the generic Gorenstein quotient of $A$ has $h$-vector $h^{''}\geq (1,2,3,4,3,2,1)$.\\\indent
Instead, if we consider the generic Gorenstein quotient of $A$ as a generic quotient of $B$, from the Main Theorem (2.9) applied to $h^{'}$ we get the sharper lower-bound $h^{''}\geq (1,3,4,4,4,3,1)$.\\
\\\indent
Since we have some substantial information on level $h$-vectors of type 2 and codimension 3, the Main Theorem (2.9) can give us some constraints on the possible level $h$-vectors of type 3 and codimension 3.\\
\\\indent
{\bf Example-Remark 2.13.} Let $h=(1,3,...,4,3)$ be the $h$-vector of a level algebra $A$. The Main Theorem (2.9) guarantees the existence of level subalgebras of $A$ of type 2 and the same socle degree having $h$-vector $h^{'}\geq (1,3,...,3,2)$.\\\indent
Hence, by $[GHMS]$'s Theorem 1.5, $h$ can be written as the sum of a level $h$-vector of the form $h^{'}=(1,3,...,4,2)$ plus $(0,0,...,0,1)$ or as a level $h$-vector of the form $h^{''}=(1,3,...,3,2)$ plus $(0,0,...,0,1,1,...,1)$.\\\indent
Again by the Main Theorem, we have that $h^{'}$ decomposes as a Gorenstein $h$-vector of codimension 3 plus $(0,0,...,0,1,1,...,1)$. Furthermore, by the Main Theorem, $h^{''}$ must decompose as the sum of a Gorenstein $h$-vector of codimension 3 and $(0,0,...,0,1)$ (see $[Za]$, Theorem 2.9 for a characterization of the level $h$-vectors of the form $(1,3,...,3,2)$).\\\indent
It follows that, in either case, $h$ can be written as a Gorenstein $h$-vector of codimension 3 plus $(0,0,...,0,1,1,...,1,2)$.\\\indent
In particular, this implies the interesting fact that level $h$-vectors of the form $h=(1,3,...,4,3)$ are {\it unimodal} (i.e. they do not increase after they start decreasing).\\\indent
Some evidence (see the tables in the Appendices of $[GHMS]$) suggests that level $h$-vectors of type 2 and codimension 3 are unimodal with their increasing part being {\it differentiable} (i.e., its first difference is an $O$-sequence). We wonder if this conjecture can be extended to level $h$-vectors of codimension 3 and arbitrary type $t$. The particular case studied in this example goes in that direction, but at this stage it seems too early to draw any conclusion. The Main Theorem (2.9) will definitely help more for $t=3$ than it does now when we have a better knowledge of the possible level $h$-vectors of type 2.\\
\\\indent
Another interesting consequence of the Main Theorem (2.9) is the following lower-bound for the next to last entry $a$ of a level $h$-vector $h=(1,r,...,a,t)$ of arbitrary type $t$: we will show in Theorem 2.14 below that, if $r\leq 7$, then $a\geq r$. It would be interesting to find a sharp lower-bound for $a$ for any given $r$ and $t$. (See $[Za]$, Theorem 3.3, ii), where we have already shown the result below for $t=2$.)\\
\\\indent
{\bf Theorem 2.14.} {\it Let $h=(1,r,...,a,t)$ be a level $h$-vector of type $t$. If $r\leq 7$, then $a\geq r$}.\\
\\\indent
{\bf Proof.} Let us proceed by induction on $t$. The Gorenstein case is obvious by symmetry, and the case $t=2$, as we just said, was shown in $[Za]$, Theorem 3.3, ii). Hence suppose that the theorem holds for type $t-1$, where $t\geq 3$. We want to prove the theorem for type $t$. Suppose, by contradiction, that, for some $r\leq 7$, there exists a level algebra $A$ having $h$-vector $h=(1,r,...,a,t)$ with $a<r$.\\\indent
We first want to show that a generic level subalgebra of $A$ of type $t-1$ and the same socle degree $e$ has $h$-vector $H^{t-1,gen}$ satisfying $H_{e-1}^{t-1,gen}\geq a-1$. Indeed, by the Main Theorem (2.9) we have $$H_{e-1}^{t-1,gen}\geq {1\over t^2-1}(r+(t^2-t-1)a)=a-{ta-r\over t^2-1},$$ which is easily seen to be greater than $a-2$ under our hypotheses. Therefore $H_{e-1}^{t-1,gen}\geq a-1$, as we desired.\\\indent
Hence, by $[GHMS]$'s Theorem 1.5, $$H_1^{t-1,gen}
\cases{\geq r-1,&if $H_{e-1}^{t-1,gen}= a-1,$\cr
	=r,&if $H_{e-1}^{t-1,gen}= a$.\cr}
$$\indent
In either case, since $r>a$, we get a level $h$-vector of type $t-1$ having the next to last entry smaller than the codimension, but this contradicts the induction hypothesis. This proves the theorem.{\ }{\ }\qed \\
\\\indent
{\bf Remark 2.15.} As  we mentioned in the Introduction, the results contained in this paper also hold when $k$ is an infinite field of characteristic $p$. In fact, if $p>e$, then no modification is required in the statements or in the proofs (since no problem arises with the coefficients of the partial derivatives of the forms of degree less than or equal to $e$). Instead, when $p\leq e$, one uses {\it divided powers} in place of inverse systems (and the action of {\it contraction} in place of that of derivation), and the results remain valid (see $[IK]$, Appendix A).\\
\\\indent
{\bf Acknowledgements.} We warmly thank Professor A. Iarrobino for sending us his preprint $[Ia]$, and for several suggestions which improved the exposition of this paper. Also, we thank Professor A. Conca for the proof of Lemma 2.7.\\\indent
Furthermore, we wish to thank the referee, whose comments have improved the presentation of this article.\\\indent
A special \lq \lq thank you" goes to our former Ph.D. Thesis supervisor, Professor A.V. Geramita, for his continuous interest and for the invaluable conversations we have had (and we continue to have) with him.

\clearpage

{\bf \huge References}\\
\\
$[BH]$ {\ } W. Bruns and J. Herzog: {\it Cohen-Macaulay rings}, Cambridge studies in advanced mathematics, No. 39, Revised edition (1998), Cambridge, U.K..\\
$[Ge]$ {\ } A.V. Geramita: {\it Inverse Systems of Fat Points: Waring's Problem, Secant Varieties and Veronese Varieties and Parametric Spaces of Gorenstein Ideals}, Queen's Papers in Pure and Applied Mathematics, No. 102, The Curves Seminar at Queen's (1996), Vol. X, 3-114.\\
$[GHMS]$ {\ } A.V. Geramita, T. Harima, J. Migliore and Y.S. Shin: {\it The Hilbert Function of a Level Algebra}, Memoirs of the Amer. Math. Soc., to appear.\\
$[Ia]$ {\ } A. Iarrobino: {\it Hilbert functions of Gorenstein algebras associated to a pencil of forms}, preprint (math.AC/0412361); Proc. of the Conference on Projective Varieties with Unexpected Properties (Siena 2004), C. Ciliberto {\it et al.} eds. (de Gruyter), to appear.\\
$[Ia2]$ {\ } A. Iarrobino: {\it Compressed Algebras: Artin algebras having given socle degrees and maximal length}, Trans. Amer. Math. Soc. 285 (1984), 337-378.\\
$[IK]$ {\ } A. Iarrobino and V. Kanev: {\it Power Sums, Gorenstein Algebras, and Determinantal Loci}, Springer Lecture Notes in Mathematics (1999), No. 1721, Springer, Heidelberg.\\
$[NW]$ {\ } N. Nisam and A. Wigderson: {\it Lower bounds on arithmetic circuits via partial derivatives}, Comput. Complexity 6 (1996/1997), No. 3, 217-234.\\
$[Za]$ {\ } F. Zanello: {\it Level algebras of type 2}, preprint (math.AC/0411228).

}

\end{document}